\documentclass[]{amsart}

\usepackage{amsmath}
\usepackage{amsthm}
\usepackage{amssymb}

\usepackage{amsfonts}

\usepackage[pdftitle={A maj-inv bijection for C_2 wr S_n},pdfauthor={Dan Bernstein},bookmarks=true]{hyperref}

\usepackage{color}

\newcommand{\Z}{\mathbb{Z}}
\newcommand{\st}{\mid} % "such that"

\providecommand{\abs}[1]{\lvert #1 \rvert}

\newcommand{\rev}{\mathbf{r}} %
\newcommand{\rtlPhi}{\overleftarrow\Phi}    %
\DeclareMathOperator{\Des}{Des} %
\DeclareMathOperator{\des}{des} %
\DeclareMathOperator{\del}{del} %
\DeclareMathOperator{\inv}{inv} %
\DeclareMathOperator{\maj}{maj} %
\DeclareMathOperator{\nrmaj}{nrmaj} %
\DeclareMathOperator{\Neg}{Neg} %
\DeclareMathOperator{\rmaj}{rmaj} %

\DeclareMathOperator{\sort}{\mathit{sort}} %

\newtheorem{thm}{Theorem}[section]
\newtheorem{prop}[thm]{Proposition}
\newtheorem{lem}[thm]{Lemma}
\newtheorem{cor}[thm]{Corollary}

\theoremstyle{definition}
\newtheorem{defn}[thm]{Definition}
\newtheorem{rem}[thm]{Remark}
\newtheorem{exmp}[thm]{Example}
\newtheorem{algo}[thm]{Algorithm}

\title{A $\maj$-$\inv$ bijection for $C_2\wr A_n$}
\author{Dan Bernstein}
\date{September 11, 2005}
\begin{document}

\bibliographystyle{alpha}

\begin{abstract}
We give a bijective proof of the MacMahon-type equidistribution over the group of
signed even permutations $C_2 \wr A_n$ that was stated in~[Bernstein. Electron. J. Combin. 11 (2004) 83]. This is
done by generalizing the bijection that was introduced in the bijective
proof of the equidistribution
over the alternating group $A_n$ in~[Bernstein and Regev. S{\'e}m. Lothar. Combin. 53 (2005) B53b].
\end{abstract}

\maketitle

\section{Introduction}

In~\cite{macmahon:indices} MacMahon proved that two \emph{permutation statistics\/}, namely
the \emph{length\/} (or \emph{inversion number\/})
and the \emph{major index\/}, are equidistributed
over the symmetric group $S_n$ for every $n>0$ (see also~\cite{macmahon:combinatory}). The question of finding a bijective proof of
this remarkable fact arose naturally. That open problem was finally solved by Foata~\cite{foata:netto},
who gave a canonical bijection on $S_n$, for each $n$, that maps one statistic to the other.
In~\cite{foSch:major}, Foata and Sch\"utzenberger proved a refinement by \emph{inverse descent classes\/} of MacMahon's theorem. The theorem has received many additional refinements and generalizations, including~\cite{carlitz:qBernoulli, carlitz:qEulerian, garsia:permutation, reiner:signed, krattenthaler:major, adin:flag, regev:wreath, regev:qstatistics, stanley:sign}.

In~\cite{adin:hyperoctahedral}, Adin, Brenti and Roichman gave an analogue of MacMahon's theorem for
the group of signed permutations $B_n = C_2 \wr S_n$. A refinement of that result by inverse descent classes appeared in~\cite{adin:equiHyperoct},
and a bijective proof was given in~\cite{foata:signedI}. These results are the ``signed'' analogues of
MacMahon's theorem, its refinement by Foata and Sch\"utzenberger and Foata's bijection, respectively.

The MacMahon equidistribution does not hold when the $S_n$ statistics are restricted to the alternating subgroups $A_n \subset S_n$. However, in~\cite{regev:alternating}, Regev and Roichman defined the $\ell_A$ (\emph{$A$-length\/}),
$\rmaj_{A_n}$ (\emph{alternating reverse major index}\/) and $\del_A$ (\emph{$A$-delent number\/})
statistics on $A_n$, and proved the following refined analogue of MacMahon's theorem:
\begin{thm}[{see~\cite[Theorem~6.1(2)]{regev:alternating}}]
For every $n>0$,
\begin{multline*}
    \sum_{w \in A_{n+1}} q^{\ell_A(w)} t^{\del_A(w)} = \sum_{w \in
    A_{n+1}} q^{\rmaj_{A_{n+1}}(w)} t^{\del_A(w)} \\
    = (1+2qt)(1+q+2q^2 t)\cdots(1+q+\dots+q^{n-2}+2q^{n-1}t)	.
\end{multline*}
\end{thm}
A bijective proof was later given in~\cite{bernstein:foataForAn} in the form of a mapping $\Psi:A_{n+1}\to A_{n+1}$ with the following properties.
\begin{thm}[{see~\cite[Theorem~5.8]{bernstein:foataForAn}}]\label{PR:Psi}
\begin{enumerate}
\item
The mapping $\Psi$ is a bijection of $A_{n+1}$ onto itself.
\item
For every $v\in A_{n+1}$, $\rmaj_{A_{n+1}}(v) = \ell_A(\Psi(v))$.
\item
For every $v\in A_{n+1}$, $\del_A(v)=\del_A(\Psi(v))$.
\end{enumerate}
\end{thm}

A ``signed'' analogue of the equidistribution over $A_n$ was given in~\cite{bernstein:macmahon} by defining the
$\ell_L$ (\emph{$L$-length\/}) and $\nrmaj_{L_n}$ (\emph{negative alternating reverse major index\/}) statistics on the
group of signed even permutations $L_n = C_2 \wr A_n \subset B_n$ and proving the following.
\begin{prop}[see~{\cite[Proposition~4.1]{bernstein:macmahon}}]\label{PR:1}
For every $B \subseteq [n+1]$
\begin{eqnarray*}
    \sum_{\{\, \pi \in
    L_{n+1} \st \Neg(\pi^{-1})\subseteq B\,\} }q^{\nrmaj_{L_{n+1}}(\pi)} = \sum_{\{\, \pi \in
    L_{n+1} \st \Neg(\pi^{-1})\subseteq B\,\} }q^{\ell_L(\pi)} \\
    = \prod_{i \in
    B}(1+q^i)\prod_{i=1}^{n-1}(1+q+\dots+q^{i-1}+2q^i) ,
\end{eqnarray*}
where $\Neg(\pi^{-1}) = \{\, -\pi(i) \st 1 \le i \le n+1,\;\pi(i)<0 \,\}$.
\end{prop}

The main result in this note is a bijective proof of Proposition~\ref{PR:1}. It is accomplished by defining a mapping $\Theta : L_{n+1} \to L_{n+1}$ for every $n>0$ and proving the following theorem.
\begin{thm}[see Theorem~\ref{TH:main}]
The mapping $\Theta$ is a bijection of $L_{n+1}$ onto itself, and for every $\pi \in L_{n+1}$, 
$\nrmaj_{L_{n+1}}(\pi) = \ell_L(\Theta(\pi))$ and $\Neg(\pi^{-1}) = \Neg(\Theta(\pi)^{-1})$.
\end{thm}

The rest of this note is organized as follows: in Section~\ref{SEC:bg} we introduce some definitions and notations and give
necessary background. In Section~\ref{SEC:decomp} we review the definition of the bijection $\Psi$ and the Main Lemma of~\cite{bernstein:macmahon}, which gives a unique decomposition of elements of $L_n$. In Section~\ref{SEC:main} we define the bijection $\Theta$ and prove the main result.

\section{Background and notation}\label{SEC:bg}

\subsection{Notation}
For an integer $a\ge 0$, let $[a]=\{1,2,\dots,a\}$ (where
$[0]=\emptyset$). Let $C_k$ be the cyclic group of order $k$, let $S_n$ be the symmetric group acting on $1,\dots,n$, and let $A_n \subset S_n$ denote the alternating group.

\subsection{The symmetric group}
Recall that $S_n$ is a Coxeter
group of type $A$, its Coxeter generators being the adjacent
transpositions $\{\,s_i\,\}_{i=1}^{n-1}$ where $s_i:=(i,i+1)$. The
defining relations are the Moore-Coxeter relations:
\[
\begin{split}
s_i^2 = 1 &\quad (1\le i \le n-1),\\
(s_i s_{i+1})^3 = 1 &\quad (1 \le i < n-1),\\
(s_i s_j)^2 = 1 &\quad (|i-j|>1).
\end{split}
\]

For every $j>0$, let
\[
    R^S_j = \{ 1,\,s_j,\, s_j s_{j-1},\,\dots,\,s_j s_{j-1} \cdots s_1 \} \subseteq S_{j+1} .
\]
Recall the following fact.
\begin{thm}[{see~\cite[pp.~61--62]{goldschmidt:characters}}]\label{THM:SCanRep}
Let $w \in S_n$. Then there exist unique elements $w_j \in R_j^S$, $1 \le j \le n-1$, such that $w = w_1 \cdots w_{n-1}$. Thus, the presentation $w=w_1\cdots w_{n-1}$ is unique. Call that presentation {\em the $S$-canonical presentation of $w$}.
\end{thm}

\subsection{The hyperoctahedral group}
\emph{The hyperoctahedral group\/} $B_n := C_2 \wr
S_n$ is the group of all bijections $\sigma$ of $\{\pm 1,\pm 2,\dots,\pm
n\}$ to itself satisfying $\sigma(-i)=-\sigma(i)$, with function composition as the group operation. It is also known as the group
of \emph{signed permutations\/}.

For $\sigma\in B_n$, we shall use \emph{window notation\/}, writing $\sigma=[\sigma_1,\dots,\sigma_n]$ to
mean that $\sigma(i)=\sigma_i$ for $i\in [n]$,
and let $\Neg(\sigma) := \{\,i\in[n] \st \sigma(i)<0 \,\}$.

$B_n$ is a Coxeter
group of type $B$, generated by $s_1,\dots,s_{n-1}$ together with
an exceptional generator $s_0:=[-1,2,3,\dots,n]$
(see~\cite[Section~8.1]{bjorner:combinatorics}). In addition to the above relations between $s_1,\dots,s_{n-1}$, we have: $s_0^2 = 1$,
$(s_0 s_1)^4 = 1$, and $s_0 s_i = s_i s_0$ for all $1<i<n$.

\subsection{The alternating group}

Let $a_i := s_1 s_{i+1}$, $1 \le i \le n-1$. Then the set $A = \{\,a_i\,\}_{i=1}^{n-1}$ generates the alternating group $A_{n+1}$. This
generating set comes from~\cite{mitsuhashi:alternating}, where it is shown that the
generators satisfy the relations
\[
\begin{split}
a_1^3 = 1,&\\
a_i^2 = 1   &\quad (1 < i \le n-1),\\
(a_i a_{i+1})^3 = 1 &\quad (1 \le i < n-1),\\
(a_i a_j)^2 = 1 &\quad (|i-j|>1)
\end{split}
\]
(see~\cite[Proposition~2.5]{mitsuhashi:alternating}).

For every $j>0$, let
\[
    R_j^A = \{1,\,a_j,\,a_j a_{j-1},\,\dots,\,a_j \cdots a_2,\,a_j \cdots a_2 a_1,\,a_j \cdots a_2 a_1^{-1}\} \subseteq A_{j+2}
\]
(for example, $R_3^A = \{1, a_3, a_3 a_2, a_3 a_2 a_1, a_3 a_2 a_1^{-1}\}$). One has the following
\begin{thm}[{see~\cite[Theorem~3.4]{regev:alternating}}]\label{THM:ACanRep}
Let $v \in A_{n+1}$. Then there exist unique elements $v_j \in R_j^A$, $1 \le j \le n-1$, such that $v = v_1 \cdots v_{n-1}$, and this presentation is unique.
Call that presentation {\em the $A$-canonical presentation of $v$}.
\end{thm}

\subsection{The group of signed even permutations}\label{SEC:L}
Our main result concerns the group $L_n:=C_2
\wr A_n$. It is the subgroup of $B_n$ of index 2 containing the
\emph{signed even permutations\/}.

For a more detailed discussion of $L_n$, see~\cite[Section~3]{bernstein:macmahon}

\subsection{$B_n$, $A_{n+1}$ and $L_{n+1}$ statistics}

Let $r=x_1 x_2\dots x_m$ be an $m$-letter word on a linearly-ordered alphabet $X$. The \emph{inversion
number\/} of $r$ is defined as \[\inv(r):=\#\{\,1\le i<j \le m \st
x_i>x_j\,\} ,\] its \emph{descent set\/} is defined as
\[
	\Des(r) := \{\,1 \le i < m \st x_i>x_{i+1}\,\} ,
\] and its \emph{descent number\/} as
\[
	\des(r) := \abs{\Des(r)} .
\]
For example, with $X=\Z$ with the usual order on the integers, if $r = 3,-4,2,1,5,-6$, then $\inv(r) = 8$, $\Des(r) = \{1,\,3,\,5\}$ and $\des(r) = 3$.

It is well known that if $w\in S_n$ then $\inv(w) = \ell_S(w)$, where $\ell_S(w)$ is the
\emph{length\/} of $w$ with respect to the Coxeter generators of $S_n$, and that
$\Des(w) = \Des_S(w):= \{\,1 \le i < n \st \ell_S(w s_i) < \ell_S(w)\,\}$, which is the descent set of $w$ in the Coxeter sense.

Define the \emph{$B$-length\/} of $\sigma \in B_n$ in
the usual way, i.e., $\ell_B(\sigma)$ is the length of $\sigma$
with respect to the Coxeter generators of $B_n$.

The $B$-length can be computed in a combinatorial way as
\[
\ell_B(\sigma) = \inv(\sigma) + \sum_{i \in \Neg(\sigma^{-1})} i
\]
(see, for example,~{\cite[Section~8.1]{bjorner:combinatorics}}).

Given $\sigma\in B_n$, the \emph{$B$-delent number\/} of $\sigma$, $\del_B(\sigma)$, is defined as the number of left-to-right minima in $\sigma$, namely
\[
\del_B(\sigma) := \#\{\,2\le j \le n \st \text{$\sigma(i)>\sigma(j)$ for all $1 \le i < j$} \,\} .
\]
For example, the left-to-right minima of $\sigma=[5,\,-1,\,2,\,-3,\,4]$ are $\{2,\,4\}$, so
$\del_B(\sigma)=2$.

The \emph{$A$-length} statistic on $A_{n+1}$ was defined in~\cite{regev:alternating} as the length of the $A$-canonical presentation. Given $v \in A_{n+1}$, $\ell_A(v)$ can be computed directly as
\begin{equation}\label{EQ:Alen}
	\ell_A(v) = \ell_S(v)-\del_S(v) = \inv(v)-\del_B(v)
\end{equation}
(see~\cite[Proposition~4.4]{regev:alternating}).

\begin{defn}[{see~\cite[Definition~3.15]{bernstein:macmahon}}]\label{DEF:Llen}
Let $\sigma \in B_n$. Define the \emph{$L$-length of $\sigma$\/} by
\[
	\ell_L(\sigma) = \ell_B(\sigma)-\del_B(\sigma) = \inv(\sigma)-\del_B(\sigma)+\sum_{i\in\Neg(\sigma^{-1})}i .
\]
\end{defn}

Given $\pi \in L_{n+1}$, let
\[
    \Des_A(\pi) := \{\, 1 \le i \le n-1 \st \ell_L(\pi a_i) \le \ell_L(\pi) \,\}	,
\]
\[
    \rmaj_{L_{n+1}}(\pi) := \sum_{i \in \Des_A(\pi)} (n-i) ,
\]
and
\[
    \nrmaj_{L_{n+1}}(\pi) := \rmaj_{L_{n+1}}(\pi) + \sum_{i\in \Neg(\pi^{-1})}i .
\]
For example, if $\pi=[5,-1,2,-3,4]$ then $\Des_A(\pi) = \{1,2\}$,
$\rmaj_{L_5}(\pi)=5$, and $\nrmaj_{L_5}(\pi)=5+1+3=9$.

\begin{rem}\label{RE:coincide}
Restricted to $A_{n+1}$, the $\rmaj_{L_{n+1}}$ statistic coincides with the $\rmaj_{A_{n+1}}$ statistic as defined in~\cite{regev:alternating} and used in Theorem~\ref{PR:Psi}.
\end{rem}

\section{The bijection $\Psi$ and the decomposition lemma}\label{SEC:decomp}
\subsection{The {F}oata bijection}

The {\em second fundamental transformation on words\/} $\Phi$ was introduced in~\cite{foata:netto} (for a full description, see~\cite[Section~10.6]{lothaire:words}). It is defined on any finite word $r=x_1 x_2 \dots x_m$ whose letters $x_1,\dots,x_m$ belong to a totally ordered alphabet. Instead of the original recursive definition, we give the algorithmic description of $\Phi$ from~\cite{foSch:major}.

\begin{algo}[$\Phi$]\label{ALGO:Phi}
Let $r=x_1 x_2 \dots x_m$ ;

1. Let $i:=1$, $r'_i := x_1$ ;

2. If $i=m$, let $\Phi(r):=r'_i$ and stop; else continue;

3. If the last letter of $r'_i$ is less than or equal to (respectively greater than) $x_{i+1}$, cut $r'_i$ after every letter less than or equal to (respectively greater than) $x_{i+1}$ ;

4. In each compartment of $r'_i$ determined by the previous cuts,
move the last letter in the compartment to the beginning of it;
let $t'_i$ be the word obtained after all those moves; put
$r'_{i+1} := t'_i \, x_{i+1}$ ; replace $i$ by $i+1$ and go to
step 2.
\end{algo}

\subsection{The covering map $f$ and its local inverses $g_u$}
Recall the $S$- and $A$-canonical presentations from Theorems~\ref{THM:SCanRep} and~\ref{THM:ACanRep}.
The following {\em covering map\/} $f$, which plays an important role in the construction of the
bijection $\Psi$, relates between $S_n$ and $A_{n+1}$ by canonical presentations.
\begin{defn}[{see~\cite[Definition~5.1]{regev:alternating}}]
Define $f:R_j^A \to R_j^S$ by
\begin{enumerate}
\item
 $f(a_j a_{j-1}\cdots a_\ell) = s_j s_{j-1}\cdots s_\ell$ if $\ell\ge
 2$, and
 \item
   $f(a_j\cdots a_1) = f(a_j\cdots a_1^{-1}) = s_j\cdots s_1$.
\end{enumerate}
Now extend $f:A_{n+1} \to S_n$ as follows:
let $v\in A_{n+1}$, $v=v_1 \cdots v_{n-1}$ its $A$-canonical presentation, then
\[
    f(v) := f(v_1)\cdots f(v_{n-1}),
\]
which is clearly the $S$-canonical presentation of $f(v)$.
\end{defn}
In other words, given $v\in A_{n+1}$ in canonical presentation $v=a_{i_1}^{\epsilon_1} a_{i_2}^{\epsilon_2} \cdots a_{i_r}^{\epsilon_r}$, we obtain $f(v)$ simply by replacing each $a$ by an $s$ (and deleting the exponents): $f(v) = s_{i_1} s_{i_2} \cdots s_{i_r}$.

The following maps serve as ``local inverses'' of $f$.
\begin{defn}
For $u \in A_{n+1}$ with $A$-canonical presentation $u= u_1 u_2 \cdots u_{n-1}$,
define $g_u:R_j^S \to R_j^A$ by
\[
    g_u(s_j s_{j-1}\cdots s_\ell) = a_j a_{j-1} \cdots a_\ell
    \quad \text{if \;$\ell\ge 2$,\; and} \quad
        g_u(s_j s_{j-1}\cdots s_1) = u_j.
\]
Now extend $g_u:S_n \to A_{n+1}$ as follows:
let $w\in S_n$, $w=w_1 \cdots w_{n-1}$ its $S$-canonical presentation, then
\[
    g_u(w) := g_u(w_1)\cdots g_u(w_{n-1}),
\]
which is clearly the $A$-canonical presentation of $g_u(w)$.
\end{defn}

\subsection{The bijection $\Psi$}

Let $w= x_1 x_2 \dots x_m$ be an $m$-letter word on some alphabet $X$. Denote
the {\em reverse\/} of $w$ by $\rev(w):=x_m x_{m-1} \dots x_1$, and let $\rtlPhi := \rev \Phi \rev$, the
{\em right-to-left Foata transformation}.

\begin{defn}\label{DEF:Psi}
Define $\Psi:A_{n+1} \to A_{n+1}$ by $\Psi(v) =
g_v(\rtlPhi(f(v)))$ .
\end{defn}

That is, the image of $v$ under $\Psi$ is obtained by applying
$\rtlPhi$ to $f(v)$ in $S_n$, then using $g_v$ as an ``inverse''
of $f$ in order to ``lift'' the result back to $A_{n+1}$.

Some of the key properties of $\Psi$ are given in Theorem~\ref{PR:Psi}.

\subsection{The decomposition lemma}
\begin{defn}
Let $r=x_1\dots x_m$ be an $m$-letter word on a linearly-ordered alphabet $X$. Define $\sort(r)$ to be the non-decreasing word with the letters of $r$.
\end{defn}
For example, with $X=\Z$ with the usual order on the integers,\\
$\sort(-4,\,2,\,3,\,-5,\,1,\,2) = -5,\,-4,\,1,\,2,\,2,\,3$.

\begin{defn}
For $\pi \in L_{n+1}$, define $s(\pi) \in L_{n+1}$ by
\[
	s(\pi) = \begin{cases}
			\sort(\pi),	&\text{if $\sum_{i \in \Neg(\pi^{-1})} i$ is even};\\
			\sort(\pi) s_1,	&\text{otherwise}.
	         \end{cases}
\]
\end{defn}

The following lemma gives a unique decomposition of every element in $L_n$ into a descent-free factor and a signless
even factor.

\begin{lem}\label{LE:main}
For every $\pi \in L_{n+1}$, the only $\sigma \in L_{n+1}$ such that $\sigma^{-1}\pi \in A_{n+1}$ and $\des_A(\sigma)=0$
is $\sigma = s(\pi)$. Moreover, $\sigma = s(\pi)$ and $u = \sigma^{-1}\pi$ satisfy $\Des_A(u)=\Des_A(\pi)$, $\inv(u)-\del_B(u) = \inv(\pi)-\del_B(\pi)$, and $\Neg(\pi^{-1}) = \Neg(\sigma^{-1})$.
\end{lem}
See~\cite[Lemma~4.6]{bernstein:macmahon} for the proof.
\begin{cor}\label{COR:uniq}
If $\sigma \in L_{n+1}$ and $\des_A(\sigma) = 0$, then for every $u \in A_{n+1}$, $s(\sigma u) = \sigma$.
\end{cor}

\section{The main result}\label{SEC:main}

\begin{defn}
Define $\Theta:L_{n+1}\to L_{n+1}$ for each $n>0$ by
\[
	\Theta(\pi) = s(\pi) \Psi( s(\pi)^{-1} \pi )	.
\]
\end{defn}

\begin{thm}\label{TH:main}
The mapping $\Theta$ is a bijection of $L_{n+1}$ onto itself, and for every $\pi \in L_{n+1}$, 
$\nrmaj_{L_{n+1}}(\pi) = \ell_L(\Theta(\pi))$ and $\Neg(\pi^{-1}) = \Neg(\Theta(\pi)^{-1})$.
\end{thm}

\begin{exmp}
As an example, let $\pi = [3,\,-6,\,-4,\,5,\,2,\,-1] \in L_6$. We have $\Des_A(\pi) = \{1,\,3,\,4\}$
and therefore $\nrmaj_{L_6}(\pi) = 4+2+1+6+4+1 = 18$. Since $\sum_{i\in\Neg(\pi^{-1})}i = 11$ is odd,
we have $\sigma:=s(\pi)=\sort(\pi)s_1=[-4,\,-6,\,-1,\,2,\,3,\,5]$ and $u:=\sigma^{-1}\pi = [5,\,2,\,1,\,6,\,4,\,3]$. One can verify that the $A$-canonical presentation of $u$ is
$u = (1)(a_2)(a_3 a_2 a_1^{-1})(a_4 a_3)$, so $f(u)=(1)(s_2)(s_3 s_2 s_1)(s_4 s_3) = [4,\,1,\,5,\,3,\,2]$. Next
we compute $\rtlPhi(f(u))$ as follows: $r:=\rev(f(u)) = [2,\,3,\,5,\,1,\,4]$. Applying Algorithm~\ref{ALGO:Phi}
to $r$ we get
\[
\begin{aligned}
r'_1 &=         2 \mid\\
r'_2 &=         2 \mid 3 \mid\\
r'_3 &=         2 \mid 3 \mid 5 \mid\\
r'_4 &=         2 \mid 3 \mid 5\;\;\;1 \mid\\
\Phi(r) = r'_5 &=         2\;\;\;3\;\;\;1\;\;\;5\;\;\;4 \quad,
\end{aligned}
\]
so $v:=\rtlPhi(f(u)) = [4,\,5,\,1,\,3,\,2]$, whose $S$-canonical presentation is\\ $v=(1)(s_2)(s_3 s_2 s_1)(s_4 s_3 s_2)$. Therefore $\Psi(u)=g_u(v) = (1)(a_2)(a_3 a_2 a_1^{-1})(a_4 a_3 a_2) = [2,\,5,\,6,\,1,\,4,\,3]$.
Finally, $\Theta(\pi) = \sigma\Psi(u) = [-6,\,3,\,5,\,-4,\,2,\,-1]$, and indeed $\ell_L(\Theta(\pi)) =
7-0+11 = 18 = \nrmaj_{L_6}(\pi)$.
\end{exmp}

\begin{proof}[Proof of Theorem~\ref{TH:main}]
The bijectivity of $\Theta$ follows from the bijectivity of $\Psi$ together with Corollary~\ref{COR:uniq}.
 
Let $\pi \in L_{n+1}$, $\sigma = s(\pi)$ and $u = \sigma^{-1}\pi$. By Definition~\ref{DEF:Llen},
\[
	\ell_L(\Theta(\pi)) = \ell_L(\sigma \Psi(u)) = \inv(\sigma\Psi(u))-\del_B(\sigma\Psi(u))+\sum_{i\in\Neg((\sigma\Psi(u))^{-1})} i .
\]
By Corollary~\ref{COR:uniq} and Lemma~\ref{LE:main},
\[
	\inv(\sigma\Psi(u))-\del_B(\sigma\Psi(u)) = \inv(\Psi(u))-\del_B(\Psi(u))
\]
and
\[
	\Neg((\sigma\Psi(u))^{-1}) = \Neg(\sigma^{-1}) = \Neg(\pi^{-1}) ,
\]
so
\[
	\ell_L(\Theta(\pi)) = \inv(\Psi(u))-\del_B(\Psi(u)) + \sum_{i \in \Neg(\pi^{-1})} i .
\]
By identity~\eqref{EQ:Alen} and Theorem~\ref{PR:Psi},
\[
	\inv(\Psi(u))-\del_B(\Psi(u)) = \ell_A(\Psi(u)) = \rmaj_{A_{n+1}}(u) = \sum_{i \in \Des_A(u)} i .
\]
Again by Lemma~\ref{LE:main}, $\Des_A(u) = \Des_A(\pi)$, whence by Remark~\ref{RE:coincide}, $\rmaj_{A_{n+1}}(u) = \rmaj_{L_{n+1}}(\pi)$.
Thus
\[
	\ell_L(\Theta(\pi)) = \rmaj_{L_{n+1}}(\pi)+\sum_{i \in \Neg(\pi^{-1})} i = \nrmaj_{L_{n+1}}(\pi) . \qedhere
\]
\end{proof}

\end{document}